\documentclass[11pt,reqno]{amsart}

\usepackage{amsmath,amssymb,amsthm}
\usepackage[margin=1in]{geometry}
\usepackage[colorlinks=true,linkcolor=blue,citecolor=blue,urlcolor=blue]{hyperref}

\theoremstyle{plain}
\newtheorem{theorem}{Theorem}

\newcommand{\Aut}{\operatorname{Aut}}
\newcommand{\AGL}{\operatorname{AGL}}

\begin{document}

\title[Perfect $1$-factorisations of $K_{18}$]
{Perfect $1$-factorisations of $K_{18}$\\ with a non-trivial automorphism group}

\author{Andrei V. Ivanov}
\address{Andrei V. Ivanov, Independent Researcher, San Jose, CA, USA.
Andrei V. Ivanov died on 16 September 2026.}

\author{Leonid M. Tertitski}
\address{Leonid M. Tertitski (corresponding author), Independent Researcher, Carlsbad, CA, USA}
\email{LeonidTertitski@gmail.com}

\dedicatory{In memory of Andrei V. Ivanov (10 November 1953 -- 16 September 2026)}



\date{October 31, 2026}

\begin{abstract}
A perfect $1$-factorisation (P1F) of $K_{2n}$ is a partition of its edges into $2n-1$
perfect matchings, any two of which form a Hamiltonian cycle. Exhaustive classification
is known only up to $K_{16}$. We present a method for the exhaustive generation, up to
isomorphism, of perfect $1$-factorisations with a non-trivial automorphism group, which
reduces the search tree far enough to make $K_{18}$ and beyond feasible. It combines two ideas. For automorphism groups of
order greater than $2$, classical results of Ihrig on the automorphisms of P1Fs are used
as tools that reject almost all prescribed cycle types without search. For the order-$2$
case, which dominates, candidates are recalculated from inherited compatibility lists,
in the spirit of Gill and Wanless, with an outer loop over two further rows that keeps
those lists small, and Ihrig's results restrict the first factor. As the numerical proof
of the method we compute all $10{,}710$ P1Fs of $K_{18}$ with a non-trivial automorphism
group; before, only six were known.
\end{abstract}

\maketitle

\section{Introduction}

A \emph{$1$-factorisation} of $K_{2n}$ is a partition of its edge set into $2n-1$
perfect matchings ($1$-factors); it is \emph{perfect} (a P1F) if the union of any two
factors is a Hamiltonian cycle \cite{MendelsohnRosa,Wallis,Rosa2019}.

\emph{What was known.} Dinitz and Garnick \cite{DinitzGarnick} classified the P1Fs of
$K_{14}$ (23 classes). Gill and Wanless \cite{GillWanless} and, independently, Meszka
\cite{Meszka2020} classified those of $K_{16}$ (3155 classes, 89 of them with a
non-trivial automorphism group); both enumerate all P1Fs first. For $K_{18}$ that is out
of reach, and only six P1Fs were known, all built algebraically: the $\AGL(1,17)$
factorisation \cite{Anderson} and five from starters
\cite{Pike,BryantMaenhautWanless}. The structural theory of automorphisms of P1Fs, due to
Ihrig \cite{Ihrig1986} and Duncan and Ihrig \cite{DuncanIhrig2000}, has been available
for decades, but it had not been turned into an exhaustive search.

\emph{What is new.} The contribution of this paper is a method: prescribe an automorphism,
generate only the P1Fs it preserves, and reduce the search tree by combining known
theorems, used as tools to reject whole cases before search, with a candidate scheme
that keeps the order-$2$ case, the bottleneck, inside the memory and time of a
workstation. None of the ingredients is new by itself; what makes the computation
possible is their combination, guided by measurements of where the search actually
spends its time (Section~\ref{sec:speedups}).

\emph{The numerical proof.} With this method we computed:

\begin{theorem}\label{thm:census}
Up to isomorphism, $K_{18}$ has exactly $10{,}710$ perfect $1$-factorisations with a
non-trivial automorphism group: $531$ with $|\Aut(\mathcal F)|>2$ and $10{,}179$ with
$|\Aut(\mathcal F)|=2$.
\end{theorem}

Six of these classes were known; the other $10{,}704$, including every class with
$|\Aut|=2$, are new. Applied without change to other orders, the method also gives
complete catalogues for $K_{20}$ with $|\Aut|\ge3$ \cite{zenodoK20,zenodoK20a3}, $K_{22}$
with $|\Aut|>6$ \cite{zenodoK22} and $K_{24}$ with $|\Aut|>3$ \cite{zenodoK24}, and it
reproduces the known results for $K_{14}$ and $K_{16}$. All catalogues and the
implementation are public \cite{zenodoK18,github}.

\section{Definitions}

Let $V=\{0,\dots,2n-1\}$. An \emph{automorphism} of a P1F $\mathcal F$ is a permutation
$\alpha$ of $V$ with $\alpha(\mathcal F)=\mathcal F$. Two P1Fs are \emph{isomorphic} if
a permutation of $V$ maps one onto the other. If $\mathcal F$ is $\alpha$-invariant, its
factors split into $\alpha$-orbits $\{M,\alpha M,\alpha^2M,\dots\}$. For an involution
$\sigma$, a factor $M$ is \emph{$\sigma$-fixed} if $\sigma M=M$, and the
\emph{diagonals} of $\sigma$ are its transposition edges $\{v,\sigma v\}$. $C(\alpha)$
is the centraliser of $\alpha$ in $S_{2n}$. \emph{Row} $r$ of a P1F is the factor
containing the edge $\{0,r\}$.

\section{The tools and the cases they reject}\label{sec:theory}

Conjugate permutations have isomorphic sets of invariant P1Fs, so one permutation per
cycle type suffices. By Cauchy's theorem every non-trivial group has an element of prime
order, and a group has an element of order $\ge3$ if and only if it has one of odd prime
order or of order~$4$; the remaining groups of order $\ge4$ are elementary abelian
$2$-groups and contain a Klein four-group $V_4$. So $|\Aut|>2$ is covered by odd primes,
order~$4$ and $V_4$, and $|\Aut|=2$ by involutions. The following known results then
decide which cases need a search at all. Each is listed with its source and with what it
rejects for $K_{18}$; self-contained proofs are in the companion paper.

\begin{enumerate}
\item[T1] \emph{An involutory automorphism fixes $0$ or $2$ vertices}
      \cite{Ihrig1986}. Rejects every involution type other than $2^9$ and $2^81^2$.
\item[T2] \emph{A fixed-point-free involution is possible only if $2n\equiv2\pmod4$}
      \cite{Ihrig1986}. Keeps $2^9$ for $K_{18}$; rejects $2^{n}$ for $K_{16}$, $K_{20}$,
      $K_{24}$.
\item[T3] \emph{An automorphism of order divisible by $4$ is possible only if
      $2n\equiv2\pmod4$} \cite{Ihrig1986,DuncanIhrig2000}. Keeps order~$4$ for $K_{18}$;
      rejects it for $K_{16}$, $K_{20}$, $K_{24}$.
\item[T4] \emph{If an automorphism of odd prime order $p$ has an odd number of
      $p$-cycles, then $p\mid2n-1$} \cite{Ihrig1986}. For $K_{18}$, rejects every type of
      order $3,5,7,11,13$ with an odd number of cycles.
\item[T5] \emph{If $\sigma$ has type $2^{n-1}1^2$, it fixes exactly one factor: its
      diagonals plus the edge joining its fixed points} \cite{Ihrig1986}. Fixes the
      first factor of the $2^81^2$ search; rejects three of its four first-factor branches.
\item[T6] \emph{First-factor dichotomy: if $n$ is odd and $\sigma$ has type $2^n$, the
      number $k$ of $\sigma$-fixed factors is $1$ or $n$; for $k=1$ the fixed factor is
      the set of diagonals, for $k=n$ each fixed factor has one diagonal}
      \cite{Ihrig1986,DuncanIhrig2000}. Rejects three of the five first-factor branches
      for $2^9$ (diagonal counts $3,5,7$).
\item[T7] \emph{For $k=n$, the fixed factors induce a P1F of $K_{n+1}$ on the quotient
      by $\sigma$} \cite[Thm.~3.13(c)]{Ihrig1986}. Fixes the labelling of the first two
      rows in the $k=9$ branch.
\end{enumerate}

Without T5 and T6 the order-$2$ case, more than $95\%$ of the census, would not be
feasible; T6 alone removes about $57\%$ of the $2^81^2$ blocks before any search.

\section{The method}\label{sec:method}

\subsection{Common framework}
\emph{Exact cover by orbits.} Under a prescribed $\alpha$, a P1F is an exact cover of
the edges by $\alpha$-orbits of factors. At each node, the candidates are the factors
through the least uncovered edge that are perfect with every committed factor and whose
orbit consists of pairwise disjoint, pairwise perfect factors; a factor is built edge by
edge and abandoned as soon as it closes a short cycle with a committed one.

\emph{Isomorph rejection.} Branches related by the stabiliser of the committed factors
in $C(\alpha)$ give isomorphic results, so one candidate per orbit is explored (orderly
generation \cite{Read,McKay}); large centralisers are handled by strong generating sets
\cite{Sims1970,Holt2005}. Every completed P1F gets a canonical form, the
lexicographically least serialisation over its Hamiltonian-cycle framings, which yields
$|\Aut|$ and makes deduplication stateless: independent runs are merged by comparing
canonical forms.

\emph{Blocks.} The search is split by its first three rows. Blocks equivalent under
$C(\alpha)$ are removed in advance, so each search is a finite list known before it
starts; blocks are independent and can be resumed. A count is claimed only when its list
is exhausted: on $K_{18}$, $4{,}944$ consecutive empty blocks were followed by new
classes.

\subsection{Automorphism group of order greater than \texorpdfstring{$2$}{2}}
Here the theorems carry the method. The orders to prescribe are the odd primes, $4$
and $V_4$ (Section~\ref{sec:theory}); T2--T4 remove most of their cycle types, and each
remaining type is searched with one representative. Prescribing an automorphism of order
$\ge3$ or a group such as $V_4$ or $S_3$ constrains each factor through its whole orbit,
so these searches are short: on $K_{18}$ all of them together take less than one hour,
and they give the $531$ classes with $|\Aut|>2$.

\subsection{Automorphism group of order \texorpdfstring{$2$}{2}}
An involution constrains the search much less, and this case is the bottleneck. T1
leaves two involution types, T5 and T6 leave three first-factor branches, and within them
the search is reduced by recalculating candidates instead of generating them.

\emph{Compatibility lists.} Gill and Wanless \cite{GillWanless} precompute the factors
compatible with the first rows and a compatibility matrix between them. We use the same
idea with one change: an outer loop over all choices of row~$4$ and row~$5$, after the
three rows fixed by the block. Only for each such choice is the list of admissible
orbits built. Every entry must then be compatible with five rows instead of three, so the
list and its matrix are far smaller and fit in memory, and the passes of the loop
together cover the whole search. Below that, each node inherits its parent's list and
keeps only the entries that are edge-disjoint from, and Hamiltonian with, the orbit just
committed, so a node's candidates are a scan of a short list rather than a new search.

\emph{Involution-specific reductions.} Whether a factor is $\sigma$-fixed or paired with
its image is decided by the partner of one vertex, so this is enforced while the factor
is generated. A factor forced to be $\sigma$-fixed is built directly on the quotient
graph of $\sigma$. In the $k=n$ branch the remaining diagonals must still fit one per
remaining fixed factor, which is checked as a look-ahead.

\subsection{Measured reductions}\label{sec:speedups}
Each was measured back-to-back against the same code without it and checked to give
identical classes on reference blocks.
\begin{enumerate}
\item[(a)] Theorems T5, T6: three of five, and three of four, branches removed without
      search.
\item[(b)] Orderly checks only near the root, where they prune; deeper, the canonical form
      deduplicates: up to $29\times$.
\item[(c)] Orbit type decided during generation: about $3\times$.
\item[(d)] Quotient-graph generation of $\sigma$-fixed factors: $1{,}434$ candidates
      became $52$ at a measured position.
\item[(e)] Diagonal look-ahead in the $k=n$ branch: about $2.2\times$.
\item[(f)] Compatibility lists with the row-$4$/row-$5$ loop: $5\times$ on $K_{18}$,
      $9\times$ on the $K_{20}$ order-$3$ search.
\end{enumerate}
Together they give about two orders of magnitude on the order-$2$ search; the $K_{18}$
census takes about $23{,}600$ thread-hours.

\section{Validation}

The method, applied unchanged to $K_{14}$ and $K_{16}$, recovers exactly the $21$ and
$89$ classes with a non-trivial automorphism group of \cite{DinitzGarnick} and
\cite{GillWanless}. On $K_{18}$ the order-$2$ search independently recovers all $179$
classes of even order $|\Aut|>2$ found by the first part. Every catalogue record was
re-checked from its matrix alone by a separate checker: every pair of factors is a
Hamiltonian cycle, every edge is covered once, $|\Aut|$ is recomputed, and no class
occurs twice.

\section{Data and code}

The $10{,}710$ classes, in canonical form with $|\Aut|$, are at
\url{https://doi.org/10.5281/zenodo.22238984} \cite{zenodoK18} (CC0). The implementation,
run scripts, block lists and a regression suite are at
\url{https://github.com/LeAndr-Graph/P1F-Census} \cite{github}.

\section*{Acknowledgements}
We thank Ian M. Wanless for making the $K_{16}$ catalogue publicly available; it was our
main test. The authors used AI assistants (Anthropic's Claude through Claude Code, and
Google Antigravity) to write and optimise parts of the software to the authors'
specification, to run and log experiments, and to draft and copy-edit parts of the text.
Every result was validated as described above, and the authors take full responsibility
for the content.

\end{document}